\input amstex
\magnification=\magstep1

\documentstyle {amsppt}
\baselineskip=24truept
\pageheight {7.5in}
\NoRunningHeads

\topmatter
\title
Reaction diffusion equations with super-linear absorption:
universal bounds, uniqueness for the Cauchy problem, boundedness
of stationary solutions
\endtitle

\author  Ross G. Pinsky \endauthor
\affil
\it Technion-Israel Institute of Technology\\
Department of Mathematics\\
Haifa, 32000, Israel\\
\tt e-mail:
pinsky\@math.technion.ac.il\endaffil

\keywords semilinear parabolic and elliptic equations,
uniqueness for the
Cauchy problem, reaction-diffusion equations, universal bounds,
stationary solutions \endkeywords
\subjclass 35K15, 35K55\endsubjclass
\thanks The  research
was supported by the Fund for the Promotion of
Research at the Technion
\endthanks
\abstract

Consider classical
solutions $u\in C^2(R^n\times(0,\infty))\cap C(R^n\times[0,\infty))$
to the parabolic reaction diffusion equation
$$
\aligned&u_t =Lu+f(x,u), \ (x,t)\in R^n\times(0,\infty);\\
&u(x,0) =g(x)\ge0, \ x\in R^n;\\
&u\ge0,\endaligned
$$
where
$$
L=\sum_{i,j=1}^na_{i,j}(x)\frac{\partial^2}{\partial x_i
\partial x_j}+\sum_{i=1}^nb_i(x)\frac\partial{\partial x_i}
$$
is a non-degenerate elliptic operator, $g\in C(R^n)$ and the
reaction  term $f$
converges to $-\infty$
at a super-linear rate as $u\to\infty$.
We give a sharp minimal growth condition on $f$, independent of $L$,
in order that there exist a universal, a priori upper bound
for all solutions to the above Cauchy problem---that is, in order that there
exist a finite function
$M(x,t)$ on $R^n\times(0,\infty)$ such that
$u(x,t)\le M(x,t)$, for all solutions to the Cauchy problem.
Assuming now in addition that $f(x,0)=0$, so that $u\equiv0$ is a solution
to the Cauchy problem,
we  show that
under a similar growth condition,
an intimate relationship exists between
two seemingly disparate phenomena---namely, uniqueness for
the  Cauchy
problem  with initial data $g=0$
and the nonexistence of unbounded, stationary solutions
to the corresponding elliptic problem.
We also give a generic condition for nonexistence of nontrivial
stationary solutions.

\endabstract
\endtopmatter

\noindent \bf 1. Introduction and statement of results.\rm\
Consider classical solutions $u\in C^2(R^n\times(0,\infty))\cap
C(R^n\times[0,\infty))$ to the parabolic reaction diffusion
equation
$$
\aligned&u_t =Lu+f(x,u), \ (x,t)\in R^n\times(0,\infty);\\
&u(x,0) =g(x)\ge0, \ x\in R^n;\\
&u\ge0,\endaligned\tag 1.1
$$
where
$$
L=\sum_{i,j=1}^na_{i,j}(x)\frac{\partial^2}{\partial x_i
\partial x_j}+\sum_{i=1}^nb_i(x)\frac\partial{\partial x_i},
$$
with $a_{i,j},b_i\in C^{\alpha}(R^n)$ and $\{a_{i,j}\}$ strictly
elliptic; that is, $\sum_{i,j=1}^na_{i,j}(x)\nu_i\nu_j>0$, for all
$x\in R^n$ and $\nu\in R^n-\{0\}$. We assume that $g\in C(R^n)$.
We require  that the reaction term $f$ be locally Lipschitz in $x$
and in $u$ and  converge to $-\infty$ at a super-linear rate as
$u\to\infty$,  for each $x\in R^n$. This latter requirement will
be made more precise below.

Our first result is a sharp minimal growth condition  on $f$,
independent of $L$, in order that there exist a universal, a
priori upper bound for all solutions to the Cauchy problem
(1.1)---that is, in order that there exist a finite function
$M(x,t)$ on $R^n\times(0,\infty)$ such that $u(x,t)\le M(x,t)$,
for all solutions to (1.1). After this result, we will always
assume that $f(x,0)=0$, so that $u\equiv0$ is a solution to (1.1).
We  show that under a growth condition similar to the above one,
an intimate relationship exists between two seemingly disparate
phenomena---namely, uniqueness for the  Cauchy problem (1.1) with
initial data $g=0$ and the nonexistence of unbounded, stationary
solutions to the corresponding elliptic problem. We also give a
generic condition for nonexistence of nontrivial stationary solutions.

For $R>0$, define
$$
F_R(u)=\sup_{|x|\le R}f(x,u).
$$
We will always assume that
$$
\sup_{u>0}F_R(u)<\infty,\ \text{for all}\ R>0.\tag F-1
$$
Theorem 1 and Example 1 below  show that the following assumption
on $F_R$ is a sharp condition for the existence of such a
universal a priori upper bound, for all solutions to (1.1). Let
$\log^{(n)}x$ denote the $n$-th iterate of $\log x$ so that
$\log^{(1)}x=\log x, \log^{(2)}x=\log \log x,$ etc.
$$
\aligned&\text{For each}\ R>0,\
\text{there exist an }\ m\ge0 \ \text{and an}\ \epsilon>0\ \text{such that}\\
&\lim_{u\to\infty}\frac{F_R(u)}
{u(\prod_{i=1}^m\log^{(i)}u)^2(\log^{(m+1)}u)^{2+\epsilon}}=-\infty,
\endaligned
\tag F-2
$$
where by convention, $\prod_{i=1}^0\log^{(i)}=1.$
\medskip

\noindent\bf Remark.\rm\ $F_R$ will satisfy (F-1) and (F-2) if, for instance,
$f(x,u)=V(x)u-\gamma(x)u^p$, for $p>1$,
or if $f$ is appropriately defined for small $u$ and satisfies
$f(x,u)=
V(x)u-\gamma(x)u(\prod_{i=1}^m\log^{(i)}u)^2(\log^{(m+1)}u)^{2+\epsilon}$,
for large $u$,
where $V(x)$ is bounded on compacts and $\gamma$ is
positive and  bounded away from 0 on
compacts.
\medskip

\noindent \bf Theorem 1.\it\ Assume that (F-1) and (F-2) hold.
Then there exists a continuous   function $M(x,t)$ on
$R^n\times(0,\infty)$ such that every solution $u$ to the Cauchy
problem (1.1) satisfies $u(x,t)\le M(x,t)$, for all $x\in R^n$ and
all $t\ge0$.
\medskip\rm

The following example shows that condition (F-2) is sharp.
\medskip

\noindent \bf Example 1.\rm\
Let  $L=\frac{d ^2}{dx^2}$
and $f(x,u)=-u\left((\log u)^2+\log u\right)$, for $u\ge1$. Then
for each $l\in R$, $u_l(x)=\exp(\exp(x+l))$ solves (1.1)
(as a stationary solution).
Since $\lim_{l\to\infty}u_l(x)=\infty$,
there is no universal a priori  upper bound for all nonnegative solutions of
(1.1) for this choice of $f$.
Alternatively, if we let
$f(x,u)=-u\left((\log u)^2(\log\log u)^2+
\log u\log\log u\right)$, for $u\ge e$, then $u_l(x)=\exp(\exp(\exp(x+l)))$
solves (1.1).  More generally,
letting $u_l(x)$ denote the $(m+1)$-th   iterate of the exponent function
with argument $x+l$, then $Lu_l+f(u_l)=0$, where the function $f$ satisfies
$f(u)<
-u(\prod_{i=1}^m\log^{(i)}u)^2$, for large $u$

\medskip

\noindent \bf Remark.\rm\
Consider the ordinary differential equation
$$
v'=f(v),\ v(0)=c\ge0,\tag1.2
$$
where $f$ is a Lipschitz function satisfying
$f(0)=0$ and $\lim_{u\to\infty}f(u)=-\infty$.
The unique solution $v_c$ to (1.2) satisfies $v_c\ge0$
and is
increasing as a function of its initial condition $c$.
It is well-known and straight forward  to show that
$\lim_{c\to\infty}v_c(t)=\infty$, if $\int^\infty\frac1{-f(u)}du=\infty$,
while $v_\infty(t)\equiv\lim_{c\to\infty}u_c(t)<\infty$, for
$t>0$, if $\int^\infty\frac1{-f(u)}du<\infty$.
Thus, if the above integral
is finite,  $v_\infty$  serves
as a universal a priori upper bound for all solutions to (1.2),
while if the above integral is infinite,
there is no  such finite function.
In particular then, for the ordinary differential
equation (1.2), a universal a priori upper bound on solutions
exists when
$f(u)=-u(\prod_{i=1}^m\log^{(i)}u)(\log^{(m+1)}u)^{1+\epsilon}$,
but not when
$f(u)=-u(\prod_{i=1}^m\log^{(i)}u)$.
Comparing this with Theorem 1 and Example 1,
one sees that \it the introduction of spatial diffusion and drift
slightly increases the minimal super-linearity threshold for
the existence of a universal a priori upper bound.
\medskip\rm\

\TagsAsMath

Define now
$$
F(u)=\sup_{x\in R^n}f(x,u),
$$
and consider  the spatially uniform versions of conditions (F-1)
and (F-2):
\TagsAsMath
$$
\sup_{u>0}F(u)<\infty; \tag F\text{-}1'
$$
$$
\lim_{u\to\infty}\frac{F(u)}
{u(\prod_{i=1}^m\log^{(i)}u)^2(\log^{(m+1)}u)^{2+\epsilon}}=-\infty,
\text{for some}\ m\ge0\ \text{and some}\ \epsilon>0.
\tag F\text{-}2'
$$
\TagsAsText
Consider also the following condition:
$$
f(x,u)=0 \ \text{and}\ F(u)\ \text{is locally Lipschitz}.\tag F-3
$$
\medskip

\noindent\bf Remark.\rm\ $F$ will satisfy (F-1$'$), (F-2$'$)
and (F-3)
if, for instance, $f$ is as in the remark following (F-2)
with $V$  bounded and $\gamma$ positive and bounded away from 0.
\medskip

\noindent The above conditions    turn out to be  critical  for
certain other important phenomena. Consider the associated
elliptic equation corresponding to stationary solutions of (1.1):
$$
\aligned&LW+f(x,W)=0,\ x\in\ R^n;\\
&W\ge0.\\
\endaligned\tag1.3
$$
We will sometimes need
one of  the following two technical conditions on $f$:
$$
\aligned& G(u)\equiv\sup_{x\in R^n}\sup_{v\ge u} (f(x,v)-f(x,v-u))\
\text{is locally Lipschitz,  is negative for large}\ u\\
& \text {and satisfies}
  \int^\infty\frac1{-G(u)}du<\infty.
\endaligned
 \tag F-4a
$$
$$
H(u)\equiv\sup_{x\in R^n}\sup_{v\ge0}(f(x,u+v)-f(x,v))\ \text{satisfies}\
 \ \text{(F-1$'$) and (F-2$'$)}.\tag F-4b
$$

\medskip
\noindent \bf Remark.\rm\ Note that if $f(x,\cdot)$ is concave
for each $x\in R^n$ and $F$ satisfies (F-1$'$), (F-2$'$)
and (F-3),
then both (F-4a) and (F-4b) hold. Indeed, by concavity, the supremum
over $v$
is attained in (F-4a)  at $v=u$ and in (F-4b) at $v=0$, giving
$G(u)=H(u)=F(u)-F(0)=F(u)$. Furthermore, the integral condition in (F-4a)
holds for any function satisfying (F-2$'$).
\medskip

\noindent\bf Theorem 2.\it\  Assume that (F-1$'$),  (F-2$'$), (F-3) and (F-4a) hold.
If the trivial solution $u=0$ is the only solution to
the Cauchy problem (1.1)
with initial data $g=0$, then all solutions $W$ to
the stationary equation (1.3)
are bounded. More specifically,
$$
W(x)\le c_0,\ \text{for all}\ x\in R^n,
$$
where $c_0$ is the largest root of the equation $G(u)=0$,
and $G$ is as in (F-4a).
In particular, if $G(u)<0$, for $u>0$, then there are no nontrivial solutions to
the stationary equation (1.3).
(If $f(x,\cdot)$ is concave for each $x$, then (F-4a) is superfluous and
$G=F$.)
\medskip\rm\

\noindent \bf Remark.\rm\
Note that if under the condition in Theorem 2, one can exhibit an  unbounded
solution to the stationary equation (1.3), then Theorem 2 guarantees the
existence of a nontrivial solution to the Cauchy
problem (1.1) with 0 initial data. Examples 2 and 3 below are applications of this.
Similarly, in the case that $G(u)<0$, for $u>0$,
if one can exhibit a nontrivial   solution to the stationary equation
(1.3), then Theorem 2 guarantees the
existence of a nontrivial solution to the Cauchy
problem (1.1) with 0 initial data.
An  application of this is given on the top  of page 9.
These examples  illustrate
the utility of Theorem 2---it
is much easier to construct appropriate solutions
to (1.3) than to construct a nontrivial solution to (1.1) with
initial data $g=0$. By Theorem 2,  with appropriate conditions
on $f$,   the existence of the latter is
guaranteed by the existence of the former.
\medskip

The next result gives  conditions for uniqueness of the Cauchy
problem (1.1) with initial data $g=0$ and also for general initial data.
Consider the following growth assumption on the coefficients
of $L$:
$$
\aligned&\sum_{i,j=1}^na_{ij}(x)\nu_i\nu_j
\le C|\nu|^2(1+|x|^2)\\
&|b(x)|\le C(1+|x|),
\endaligned\tag L-1
$$
for some $C>0$.

\medskip

\noindent \bf Theorem 3.\it\ If (F-1$'$), (F-2$'$) (F-3) and (L-1)  hold, then
the trivial solution $u\equiv0$ is the only solution to
the Cauchy problem (1.1) with
initial data $g=0$. If in addition, (F-4b) holds, then
there is a unique solution to the Cauchy problem (1.1) for each $g\in C(R^n)$.
\medskip\rm\

As immediate corollaries to Theorems 2 and 3, we obtain the following theorems.

\noindent \bf Theorem 4.\it\ If (F-1$'$), (F-2$'$), (F-3), (F-4a) and (L-1) hold, then
all  solutions to the stationary equation (1.3) are bounded.

\medskip

\noindent \bf Theorem 5.\it\ Assume that
 (F-1$'$), (F-2$'$), (F-3), (F-4a) and (L-1) hold.
Assume in addition that the function $G$ from condition (F-4a)
satisfies
$G(u)<0$,
for $u>0$
(which will occur in particular if $f(x,\cdot)$
is concave for each $x$ and $F(u)<0$, for $u>0$).
Then
there are no nontrivial solutions to the stationary equation (1.3).

\medskip\rm

We elaborate now on Theorems 3 and 4 and then on Theorem 5.
We begin by providing two examples which demonstrate that
condition (L-1) is sharp for  both
Theorem 3 and Theorem 4.

\medskip\rm

\noindent \bf Example 2.\rm\
When $L=(1+x^2)^{1+\epsilon}\frac{d^2}{dx^2}$ and
$f(x,u)=-2u^{1+\epsilon}$, for  some $\epsilon>0$, then (1.3) possesses
the unbounded solution    $u(x)=1+x^2$. By Theorem 2, it then follows that there
exists a nontrivial  solution to (1.1) with initial data $g=0$.

\medskip

\noindent\bf Example 3.\rm\
When $L=\frac{d^2}{dx^2}+(1+x^2)^{\frac12+\epsilon}$sgn$(x)\frac d{d x}$
 and $f(x,u)=-2-2(1+u)^{\frac12+\epsilon}u^\frac12$, for some $\epsilon>0$,
then
(1.3) possesses the unbounded solution $u(x)=x^2$.
By Theorem 2, it then follows that there exists
a nontrivial solution  to (1.1) with initial data $g=0$.

\medskip

Theorems 2-4 and Examples 2-3 suggest that
under conditions (F-1$'$), (F-2$'$), (F-3) and the technical condition (F-4a),
which is always satisfied if $f(x,\cdot)$ is concave,
there may well be an   equivalence between
uniqueness for the Cauchy problem (1.1) with initial data $g=0$
and nonexistence of unbounded
solutions to the stationary equation (1.3).
We leave this as an open problem.
\medskip

\noindent \bf Remark.\rm\ We emphasize that in the context of this paper,
uniqueness for the Cauchy problem (1.1) means uniqueness with regard to \it all
\rm\ classical, nonnegative solutions. If one works only with, say, mild
solutions, then the situation can be quite different.
For example, there is a unique mild solution
for $u_t=\Delta u-\gamma(x)u^p$, for $p>1$
and bounded $\gamma\ge0$ \cite{8}; yet
if $\gamma$ decays sufficiently rapidly, uniqueness fails in the sense
of all classical solutions (for details,
see the next to the last paragraph in this section).

\medskip
Turning to  Theorem 5, it follows in particular
that if $f(x,\cdot)$ is concave for each $x$, $F(u)<0$, for $u>0$,
$F$ satisfies conditions  (F-2$'$) and (F-3),
and the operator $L$ satisfies
condition (L-1), then there are no nontrivial solutions to (1.3).
Generic results such as this regarding existence/nonexistence
of nontrivial solutions to (1.3) seem to be rare in the literature.
Indeed,  the question of existence/nonexistence  is  delicate
and can hinge greatly on the particular form of $L$ and $f$.
One generic result in the literature concerns the case that
$f(\cdot,u)$ and the coefficients of $L$ are  periodic in $x$.
Assume that for some $M_0>0$, $f(x,u)\le0$, for all $x\in R^n$ and all
$u\ge M_0$.
Let $\lambda_0$ denote the principal eigenvalue for the operator
$L+\frac{\partial f}{\partial u}(x,0)$ with periodic boundary conditions.
If $\lambda_0>0$,
then (1.3)  possesses a nontrivial periodic solution, while
if $\lambda_0\le 0$ and $\frac{f(x,u)}u$ is decreasing in $u$ for each $x\in R^n$,
then
(1.3) does not  possess a nontrivial \it bounded\rm\ solution
\cite{1}.

Consider now the  well-studied case  $L=\alpha(x)\Delta$
and $f(x,u)=-u^p$, with $p>1$. If $n\ge2$, then
(1.3) possesses a nontrivial solution if $\lim_{|x|\to\infty}
\frac{\alpha(x)}{(1+|x|)^{2+\epsilon}}>0$, for some $\epsilon>0$,
and does not possess a nontrivial solution if
$\lim_{|x|\to\infty}
\frac{\alpha(x)}{(1+|x|)^2}<\infty$.
For $n=1$, the same result holds with the exponent 2 replaced by
$1+p$. For $n\ge3$, this result goes back to
\cite{5}  and \cite{7}, and it is shown in \cite{5}
that in the case of existence there are in fact an infinite number
of bounded solutions.
The $n$-dimensional analog of Example 2 above shows
that there is  also an  unbounded solution.
For $n=1,2$, the above result  were proven in \cite{2}
and later apppeared with a different proof in \cite{3}
(which also re-derives the result for $n\ge3$).
Note that by Theorem 5, nonexistence of nontrivial solutions to (1.3)
continues to hold for $\alpha$ in the above nonexistence range
when the nonlinearity $-u^p$  is replaced by
$f(x,u)=-u(\log (u+1))^{2+\epsilon}$ or
$f(x,u)=-u(\log(u+e))^2(\log\log(u+e))^{2+\epsilon}$, etc.,
for some $\epsilon>0$.
Also, note that
when $\alpha$ is in the above \it existence\rm\ range,
then  by Theorem 2,
there is a nontrivial solution to the Cauchy problem
$u_t=\alpha \Delta u-u^p$ with  initial condition $g=0$.

An open problem was mentioned after Example 3.
We now discuss some more  open problems suggested by the
above results and  make some informal conjectures.
Example 1 above shows that condition (F-2$'$)
is sharp for Theorem 4.
We don't believe that condition (F-2$'$) is sharp for Theorem 3.
That is, we don't believe that the intimate connection between
uniqueness for the Cauchy problem (1.1) with initial condition $g=0$
 and nonexistence of unbounded
solutions to (1.3) continues to hold when condition (F-2$'$) is not
in effect. Indeed, considering that uniqueness holds for positive
solutions to the linear
Cauchy
problem $u_t=\Delta u-u$ and, by Theorem 3, also for the  Cauchy
problem $u_t=\Delta u+f(u)$,
when $f$ approaches $-\infty$ sufficiently fast so as to satisfy (F-2$'$),
it seems likely that uniqueness also holds
for $u_t=\Delta u+f(u)$ when $f$ approaches $-\infty$
at a super-linear rate that does not satisfy (F-2$'$).
We leave this as an open problem.
If one replaces $\Delta$ with a general operator
$L$, then the above heuristics become  more uncertain.
Indeed, for the linear equation $u_t=Lu-u$, uniqueness is known
to hold
when $b$ satisfies the condition in (L-1) and when $a$ satisfies a two-sided bound
of the form
$c
(1+|x|^\gamma)|\nu|^2\le
\sum_{i,j=1}^na_{ij}(x)\nu_i\nu_j
\le C(1+|x|^\gamma)|\nu|^2$,
for some $\gamma\in[0,2]$.
It is not known whether uniqueness holds for the linear problem under condition
(L-1). (For more about
 uniqueness in the linear case, see \cite{4} and references therein.)

We also note that condition (F-2) in place of (F-2$'$) is not sufficient
to insure uniqueness for the Cauchy problem (1.1).
An example is given in the paragraph after next.

We believe that Theorem 2, and consequently Theorem 4,
hold without the technical condition (F-4a). Similarly, Theorem 3
for nonzero initial data probably holds  without the technical
condition (F-4b). These are also open problems.

\medskip

We wish to
emphasize an important point
with regard to the connection between non-existence of unbounded
solutions to (1.3) and uniqueness for the Cauchy problem (1.1) with
initial condition $g=0$. If one takes a scalar function $\alpha(x)>0$
and replaces $L$ and $f$ by $\alpha L$ and $\alpha f$ respectively,
then of course (1.3) remains unchanged. However,
making the same change in the parabolic equation (1.1)
can affect the question of uniqueness.
Indeed, in \cite{3}, it was shown that
if $L=\Delta$ and $f(x,u)=\frac1{\alpha}u^p$ with $p>1$, then
there exists a nontrivial solution to (1.1) with initial
data $g=0$ if $\alpha(x)\ge C\exp(|x|^{2+\epsilon})$,
for some $\epsilon>0$ and $C>0$,
and there doesn't exist such a solution if
$\alpha(x)\le C\exp(|x|^2)$,
for some $C>0$.
(This is the example alluded to in the remark after Example 3
and in the paragraph before last.)
On the other hand, for $L=\alpha\Delta$ and
$f(x,u)=u^p$, it
follows from Theorem 3 and from the final sentence in the  paragraph
on the top of  page
9 that
the  existence of a nontrivial solution to (1.1) with initial data $g=0$
depends on whether $\alpha(x)\ge C(1+|x|)^{l+\epsilon}$
or $\alpha(x)\le C(1+|x|)^l$, where $l=2$, if $n\ge2$, and $l=1+p$, if $n=1$.
What allows for the connection between uniqueness for the Cauchy problem (1.1)
with initial condition $g=0$ and
nonexistence of solutions
to the stationary equation  (1.3)
is the assumption (F-2$'$) on $F(u)=\sup_{x\in R^n}f(x,u)$.
In the above example, this assumption requires one
to consider
$L=\alpha\Delta$ and $f(x,u)=u^p$, rather than
$L=\Delta$ and $f(x,u)=\frac1{\alpha}u^p$.

Theorems 1 and 3 are proved by constructing appropriate super
solutions,  the proof of Theorem 3 being the much more delicate one.
The proofs are given in sections two and three respectively.
The proof of Theorem 2 uses a gamut of techniques and is given
in section four.

\medskip

\noindent \bf 2. Proof of Theorem 1.
\rm\ We  begin with a standard maximum principle.

\noindent \bf Proposition 1. \it\
Let $D\subset R^n$ be a bounded domain and let
$0\le u_1, u_2\in C^{2,1}(D\times(0,\infty))\cap C(\bar D\times
[0,\infty))$ satisfy
$$
 Lu_1+f(x,u_1)-\frac{\partial u_1}{\partial t}
\le  Lu_2+f(x,u_2)-\frac{\partial u_2}{\partial t},\
\text{for}\  (x,t)\in D\times (0,\infty),
$$
$$
u_1(x,t)\ge u_2(x,t), \ \text{for}\ (x,t)\in
\partial D\times(0,\infty)
$$
and
$$
u_1(x,0)\ge u_2(x,0),\ \text{ for}\  x\in D.
$$
Then $u_1\ge u_2$ in $D\times(0,\infty)$.
\medskip\rm

\noindent \bf Proof.\rm\
Let $W=u_1-u_2$ and define
$V(x,t)=\frac{f(x,u_1(x,t))-f(x,u_2(x,t))}{W(x,t)}$,
if $W(x,t)\neq0$, and $V(x,t)=0$ otherwise.
Since $f$ is locally Lipschitz in $u$, $V$ is bounded in $D\times[0,T]$,
for any $T>0$.
We have
$LW+VW-\frac{\partial W}{\partial t}\le0$ in
$D\times(0,\infty)$,
$W(x,0)\ge0$ in $D$, and $W(x,t)\ge0$ on
$\partial D\times(0,\infty)$. Thus, by the standard linear
parabolic maximum principle, $u_1\ge u_2$.
\hfill
$\square$
\medskip

We record the following result, mentioned in the
remark after Example 1.
\medskip

\noindent \bf Lemma 1.\it\ Let $G(u)$ be Lipschitz
and satisfy $\lim_{u\to\infty}G(u)=-\infty$.
For $c\ge0$, let $v_c(t)$ denote the solution
to
$$
\aligned&
v'=G(v), \ t>0;\\
&v(0)=c.
\endaligned\tag2.1
$$

\noindent i. If  $\int^\infty\frac1{-G(u)}du<\infty$,
then $v_\infty(t)\equiv\lim_{c\to\infty}v_c(t)<\infty$,
for all $t>0$,  and $v_\infty$ solves (2.1) with $c=\infty$.

\noindent ii. If $\int^\infty\frac1{-G(u)}du=\infty$,
then $\lim_{c\to\infty}v_c(t)=\infty$, for all $t\ge0$.

\medskip\rm

\noindent \bf Proof.\rm\ We omit the straight forward proof of this
standard result.\hfill $\square$
\medskip

We now give the proof of Theorem 1.
It suffices to show that for some $T_0>0$ and each $R>0$, there exists
a continuous function $M_R(x,t)$ on $\{|x|<R\}\times(0,T_0]$
such that
every solution $u$ to (1.1) satisfies
$u(x,t)\le M_R(x,t)<\infty$, for $|x|<R$ and $t\in (0,T_0]$.
The reason it is enough to  consider only
$t\in(0,T_0]$ is that if $u(x,t)$ is a solution
to (1.1), then $u(x,T_0+t)$ is a solution to (1.1)
with the initial condition  $g(\cdot)$ replaced by $u(\cdot, T_0)$.

We will assume that $F_R$ satisfies (F-2) with $m=0$.
At the end of the proof, we  describe the simple change needed in the
case that $m\ge1$.
In particular then, there exists
an $\epsilon>0$ and a $u_0>1$ such that
$$
F_R(u)\le -u(\log u)^{2+\epsilon}\equiv Q(u),\ \text{for}\ u\ge u_0.\tag2.2
$$
Since $\int^\infty\frac1{-Q(u)}du<\infty$, it follows  from
Lemma 1 that there exists a $T_0>0$ and  a function $v_\infty(t)$ satisfying
$$
\aligned&v_\infty'=Q(v_\infty),\ t\in(0,T_0];\\
&v_\infty(0)=\infty;\\
&v_\infty(t)>1, \ t\in(0,T_0].\endaligned\tag2.3
$$
Define
$$
\phi_R(x)=\exp((R^2-|x|^2)^{-l}),
$$
with $l$ satisfying $l\epsilon>2$.
Finally, choose $K$ so that $\exp(K)>u_0$ and define
$$
M_R(x,t)=\exp((K(t+1))\phi_R(x)+v_\infty(t).
$$

Since $M_R(x,t)>u_0$, it follows from (2.2)
that $f(x,M_R)\le F_R(M_R)\le Q(M_R)$.
Since $Q(u)$ is concave  for $u\ge1$
and $Q(1)=0$, it follows from the mean value theorem that
$Q(b+a)-Q(b)< Q(a)$ for $1\le a\le b$.
Thus, since
$\exp((K(t+1))\phi_R(x),v_\infty(t)>1$,
we have
$Q(M_R(x,t))< Q(\exp(K(t+1))\phi_R(x))+Q(v_\infty(t))$.
Using these facts along with (2.3), we obtain
$$
\aligned&LM_R+f(x,M_R)-(M_R)_t\le
\exp(K(t+1))L\phi_R+Q(\exp(K(t+1))\phi_R)+Q(v_\infty)\\
&-K\exp(K(t+1))\phi_R-v_\infty'=\\
&\exp(K(t+1))L\phi_R+Q(\exp(K(t+1))\phi_R)-K\exp(K(t+1))\phi_R<\\
& \exp(K(t+1))\left(L\phi_R-(R^2-|x|^2)^{-(2+\epsilon)l}\phi_R-
K\phi_R\right),\ \text{for}\ |x|<R\ \text{and}\ t\in(0,T_0].
\endaligned
\tag2.4
$$
We  have
$$
\aligned&\frac{L\phi_R(x)}{\phi_R(x)}=
\left(4l^2(R^2-|x|^2)^{-2l-2}
+4l(l+1)(R^2-|x|^2)^{-l-2}\right)\sum_{i,j=1}^na_{i,j}(x)x_ix_j\\
&+2nl(R^2-|x|^2)^{-l-1}
+2l(R^2-|x|^2)^{-l-1}\sum_{i=1}^nx_ib_i(x).
\endaligned\tag2.5
$$
The right hand side of (2.5) is bounded for $|x|$ in any ball of radius
less than $R$. Furthermore,
on the right hand
side of (2.5),
the dominating term  as $|x|\to R$ is
$4l^2(R^2-|x|^2)^{-2l-2}$.
Thus, since $l\epsilon>2$,
it follows that the right hand side of (2.4) is negative if $K$
is chosen sufficiently large.
Using this with Proposition 1 and the fact that $M_R(x,0)=\infty$
and $M_R(x,t)=\infty$, for $|x|=R$, we conclude that
any solution $u$ to (1.1) satisfies
$u(x,t)\le M_R(x,t)$, for $|x|<R$ and $t\in(0,T_0]$.
This completes the proof of the theorem under the assumption
that $m=0$ in (F-2).

When $m>0$ one simply replaces the test function
$\phi_R(x)$ as above by $\phi_R(x)=\exp^{(m+1)}((R^2-|x|^2)^{-l})$,
where $\exp^{(j)}$ denotes the $j$-th iterate of the exponential function.
Everything goes through in a similar fashion.

\hfill $\square$

\medskip

\noindent \bf 3. Proof of Theorem 3.\rm\ By assumption,  $F(u)=\sup_{x\in  R^n}f(x,u)$ satisfies (F-2$'$).
As we did in the proof of Theorem 1, we will assume that  $m=0$ in (F-2$'$).
At the appropriate point in   the proof,
we  describe the simple change needed in the
case that $m\ge1$.

We first consider the case with initial condition  $g=0$.
By conditions (F-2$'$) and (F-3),
it follows that there exist
$C_0,\epsilon>0$ and $M_0>1$ such that
$$
\aligned&f(x,u)\le F(u)\le C_0u,\ \text{for}\ u\le M_0;\\
&f(x,u)\le F(u)\le-u(\log u)^{2+\epsilon},\ \text{for}\  u\ge M_0.
\endaligned\tag3.1
$$
Fix    $R>1$ and $T\in(0,\infty)$. Define
$$
\phi_R(x)=\exp((\frac{1+|x|^2}{R^2-|x|^2})^l),
$$
with $l$ satisfying $l\epsilon>2$, and define
$$
\psi_R(x,t)=(\phi_R(x)-1)\exp(K(t+1)),
$$
with $K>0$.
A direct calculation reveals that
$$
L\psi_R=
\exp(K(t+1))\phi_R(x)\left[W_1+W_2+W_3+W_4+W_5\right],\tag3.2
$$
where
$$
\aligned&W_1=4l^2(1+|x|^2)^{2l-2}(R^2-|x|^2)^{-2l-2}(R^2+1)^2
\sum_{i,j=1}^na_{i,j}(x)x_ix_j;\\
&
W_2=4l(l-1)(1+|x|^2)^{l-2}(R^2-|x|^2)^{-l-2}(R^2+1)^2\sum_{i,j=1}^na_{i,j}(x)x_ix_j;\\
&
W_3=4l(1+|x|^2)^{l-1}(R^2-|x|^2)^{-l-2}(R^2+1)\sum_{i,j=1}^na_{i,j}(x)x_ix_j;\\
&
W_4=2nl(1+|x|^2)^{l-1}(R^2-|x|^2)^{-l-1}(R^2+1);\\
&
W_5=2l(1+|x|^2)^{l-1}(R^2-|x|^2)^{-l-1}(R^2+1)\sum_{i=1}^nx_ib_i.\endaligned
$$
We also have
$$
\frac{\partial \psi_R}{\partial t}=K\psi_R(x).\tag3.3
$$

We claim that for $K$ sufficiently large and independent of $R$
(but not independent of $T$ in (3.4-b))
$$
\aligned&
\exp(K(t+1))\phi_R(x)W_i-\frac15(K-C_0)\psi_R(x,t)\le0,\\
&\text{if}\ \psi_R(x,t)\le M_0, \ \text{for}\ |x|<R\ \text{and}\ t\in (0,T];
\endaligned
\tag3.4-a
$$
$$
\aligned&
\exp(K(t+1))\phi_R(x)W_i-
\frac15K\psi_R(x,t)-\frac15\psi_R(x,t)(\log\psi_R(x,t))^{2+\epsilon}\le0,\\
& \text{if}\ \psi_R(x,t)\ge M_0,\
\text{for}\ |x|<R\ \text{and}\ t\in (0,T],\endaligned\tag3.4-b
$$
for $ i=1,2,3,4,5$

From (3.1)-(3.4), it follows that for sufficiently large $K$, independent of $R$,
$$
L\psi_R-\frac{\partial \psi_R}{\partial t}+f(x,\psi_R)\le0, \ \text{for}\
|x|<R\ \text{and}\ t\in(0,T].
\tag3.5
$$
Since $\psi_R(x,0)\ge0$ and $\lim_{|x|\to R}\psi_R(x,t)=\infty$,
it  follows from (3.5) and the maximum principle in Proposition 1
that
any solution $u$ to (1.1) with  initial condition $g=0$
must satisfy the bound
$$
u(x,t)\le
(\exp((\frac{1+|x|^2}{R^2-|x|^2})^l)-1)\exp(K(t+1)),
\ \text{for}\ |x|<R \ \text{and}\ t\in(0,T].\tag3.6
$$
Since $K$ doesn't depend on $R$, letting $R\to\infty$ in (3.6)  gives
$u(x,t)\equiv0$ for $x\in R^n$ and $t\in(0,T]$.
Now letting $T\to\infty$ gives $u(x,t)\equiv0$ in $R^n\times(0,\infty)$,
completing the proof.

When $m>0$ one  replaces the test function
$\phi_R(x)$ as above by  $\phi_R(x)=$
\linebreak $\exp^{(m)}((\frac{1+|x|^2}{R^2-|x|^2})^l)$,
where $\exp^{(j)}$ denotes the $j$-th iterate of the exponential function.
The resulting calculations are similar to the present case.

It thus   remains  to prove (3.4) for $K$ independent of $R$.
We will prove (3.4) for $W_1$.
The proofs for $W_i, i\ge2$, are similar.
Consider first (3.4-a).
We will always assume that $K\ge C_0$.
Recall the definitions of $\phi_R$ and
$\psi_R$.  If $\psi_R(x,t)\le M_0$, then
a fortiori $\phi_R(x)\le M_0+1$ and
$(\frac{1+|x|^2}{R^2-|x|^2})^l\le \log( M_0+1)\equiv L_0^l$.
Also, we have $\psi_R(x,t)\ge\phi_R(x)-1\ge (\frac{1+|x|^2}{R^2-|x|^2})^l$.
In light of these observations, it follows that (3.4-a) will
hold if
$$
(M_0+1)W_1-\frac15(K-C_0)(\frac{1+|x|^2}{R^2-|x|^2})^l\le0,
\ \text{whenever}\
\frac{1+|x|^2}{R^2-|x|^2}\le L_0.
$$
Or equivalently,
if
$$
K\ge C_0+5(\frac{R^2-|x|^2}{1+|x|^2})^l(M_0+1)W_1,
\ \text{whenever}\
\frac{1+|x|^2}{R^2-|x|^2}\le L_0.
\tag3.7
$$
Thus, we must show that the right hand side of (3.7) is bounded in $R$ and $x$
under the constraint
$\frac{1+|x|^2}{R^2-|x|^2}\le L_0$.
Substituting for $W_1$ in the right hand side of (3.7) and using the
assumption that $\sum_{i,j=1}^na_{i,j}(x)x_ix_j\le C(1+|x|^2)|x|^2$,
one finds that it is enough to show that
$(\frac{1+|x|^2}{R^2-|x|^2})^{l-1}\frac{(R^2+1)^2|x|^2}{(R^2-|x|^2)^3}$
is bounded in $R$ and $x$ under the above constraint.
Since
$(\frac{1+|x|^2}{R^2-|x|^2})^{l-1}$ is trivially bounded under the constraint,
it remains only to consider
$\frac{(R^2+1)^2|x|^2}{(R^2-|x|^2)^3}$.
The  constraint above  is equivalent to the constraint
$|x|^2\le \frac{L_0R^2-1}{L_0+1}$. From this  it is clear that
under the constraint,
$\frac{(R^2+1)^2|x|^2}{(R^2-|x|^2)^3}$ is bounded in $R$ and $x$.

We now turn to (3.4-b). The constraint
$\psi_R(x,t)\ge M_0$ along with the condition $t\le T$ guarantee the
existence of a $c_0\in(0,1)$ such that $\phi_R(x)-1\ge c_0\phi_R(x)$.
Note that $c_0$ depends on $T$, but not on $R$.
Thus, under the constraint, we have $\psi(x,t)= (\phi_R(x)-1)\exp(K(t+1))
\ge c_0\phi_R(x)\exp(K(t+1))$.
Therefore (3.4-b) will hold if we show that $K$ can be picked independent of $R$
and such that
$W_1-\frac15c_0K-\frac15c_0\left(\log\phi_R+\log c_0+K(t+1)\right)
^{2+\epsilon}\le0$
holds under the constraint.
We will always assume that $K\ge-\log c_0$.
Thus it suffices to show that
$W_1-\frac15c_0(\log\phi_R)
^{2+\epsilon}$
is bounded from above under the constraint, independent of $R$.
Substituting for $\phi_R$ and $W_1$ and using
the assumption $\sum_{i,j=1}^na_{i,j}(x)x_ix_j\le C(1+|x|^2)|x|^2$,
it is sufficient to show that
$4l^2C(1+|x|^2)^{2l-1}(R^2-|x|^2)^{-2l-2}(R^2+1)^2|x|^2
-\frac15c_0\left(\frac{1+|x|^2}{R^2-|x|^2}\right)^{(2+\epsilon) l}$ is bounded
from above under the constraint,
or equivalently, that
$$
\aligned&\left(\frac{1+|x|^2}{R^2-|x|^2}\right)^{(2+\epsilon) l}
\left(4l^2C\frac{(R^2+1)^2|x|^2}{(1+|x|^2)^3}
\left(\frac{R^2-|x|^2}{1+|x|^2}\right)^{l\epsilon -2}-\frac{c_0}5\right)\\
& \text{is bounded from above under the constraint}.\endaligned
\tag3.8
$$
We may assume that
$4l^2C\frac{(R^2+1)^2|x|^2}{(1+|x|^2)^3}
\left(\frac{R^2-|x|^2}{1+|x|^2}\right)^{l\epsilon -2}\ge \frac{c_0}5$,
since otherwise it is clear that
(3.8) holds.
From this inequality and the assumption that $l\epsilon>2$, it follows that
$$
\frac{1+|x|^2}{R^2-|x|^2}\le
\left(\frac{20l^2C}{c_0}\frac{(R^2+1)^2|x|^2}{(1+|x|^2)^3}\right)^{\frac1{l\epsilon-2}}.
\tag3.9
$$
Furthermore, the constraint $\psi_R\ge M_0$ guarantees that
$$
\frac{1+|x|^2}{R^2-|x|^2}\ge(\log(1+M_0\exp(-K(T+1))))^\frac1l\equiv \gamma_0>0,
\tag3.10
$$
which can be written in the form
$$
|x|^2\ge\frac{\gamma_0R^2-1}{\gamma_0+1}.\tag3.11
$$
If $|x|$ satisfies (3.11), then the right hand side of (3.9) is bounded.
Therefore, in (3.8), the terms
$(\frac{1+|x|^2}{R^2-|x|^2})^{(2+\epsilon) l}$
and $\frac{(R^2+1)^2|x|^2}{(1+|x|^2)^3}$
are bounded. And by (3.10), the term
$(\frac{R^2-|x|^2}{1+|x|^2})^{l\epsilon -2}$
is also bounded. This completes the proof of (3.8).

We now turn to the case that the initial condition $g$ is not equal to
0. We assume now in addition that condition
(F-4b) is in effect. Fix $R>1$ and $T\in(0,\infty)$.
Let $\psi_R(x,t)$ be as in part (i), but
corresponding to the function $H$ appearing in condition (F-4b), rather
than corresponding to the function $F$ as in part (i).

In \cite{3}, for the case $f(x,u)=V(x)u-\gamma(x)u^p$,
we showed that there exists a minimal solution $u_g$ to (1.1); that
is, a solution $u_g$ with the property that $u_g(x,t)\le u(x,t)$, for
any solution $u$ to (1.1) with initial data $g$.
In fact, the proofs there go through   for general locally
Lipschitz continuous $f$
as long as a universal a priori upper bound exists. Thus, in light of
Theorem 1,
there exists a minimal solution $u_g$. (In fact, $u_g$ is obtained
by taking the solution of (4.1) below
and letting $m\to\infty$.)

Now define $\hat\psi_R(x,t)=\psi_R(x,t)+u_g$. Then
$$
\aligned&L\hat\psi_R+f(x,\hat\psi_R)-\frac{\partial \hat\psi_R}{\partial t}=
(L\psi_R+H(\psi_R)-\frac{\partial \psi_R}{\partial t})\\
&+
(Lu_g+f(x,u_g)-\frac
{\partial u_g}{\partial t})+(f(x,\psi_R+u_g)-f(x,u_g)-H(\psi_R)).
\endaligned\tag3.12
$$
The first of the three terms  on the right
hand side of (3.12) is non-positive by the construction
in part (i), the second term is non-positive because $u_g$ is
a solution to (1.1), and the third term is non-positive by
the definition of $H$ in (F-4b).
The argument used above for part (i)
in the paragraph in which (3.5) appears then shows
that any solution $u$ to (1.1) must satisfy
$u(x,t)\le u_g(x,t)+\psi_R(x,t)$, for $|x|<R$ and $t\in(0,T]$.
Letting $R\to\infty$ and then  $T\to\infty$ as before
shows that $u=u_g$.
\hfill $\square$

\medskip

\noindent \bf 4. Proof of Theorem 2.\rm\
We need to utilize certain constructions that were carried out
in \cite{3, section 2} for the case that $f(x,u)=V(x)-\gamma(x)u^p$.
These constructions are based on
results in \cite{6}, and
hold with the same proofs for  general locally
Lipschitz continuous $f$
as long as a universal a priori upper bound exists. Thus, in light of
Theorem 1,
they hold for $f$ satisfying (F-1) and (F-2).

Let $B_m\subset R^n$ denote the open ball of radius $m$ centered at the origin.
There exists a solution
$u\in C^{2,1}
(B_m\times(0,\infty))\cap C(B_m\times[0,\infty))\cap C(\bar B_m\times(0,\infty))$
to the equation
$$
\aligned&u_t=Lu+f(x,u), \ (x,t)\in B_m\times(0,\infty);\\
&u(x,0)=g(x), \ x\in B_m;\\
&u(x,t)=0,\ (x,t)\in \partial B_m\times(0,\infty),
\endaligned\tag4.1
$$
for any  $0\le g\in C(\bar B_m)$.
(See the beginning of the proof of Theorem 1 in \cite{3}, where
the above construction is first made in the case that $g$ is compactly
supported in $B_m$, and then extended to the case
that $g\in C(\bar B_m)$.)

Now let $W$ be an arbitrary solution to (1.3).
For $m>0$ and a positive integer $k$, let $\psi_{m,k} \in
C^\infty(R^n)$    satisfy
$$
\gather \psi_{m,k}(x)=0, \ |x|\le m \
\text{and}\ |x|>2m+1\\ \psi_{m,k}(x)=k,\ m+1\le|x|\le
2m\\ 0\le \psi_{m,k}\le k.\endgather
$$
There exists  a nonnegative solution
$U_{m,k}\in C^{2,1} (B_{2m}\times(0,\infty))\cap C(\bar
B_{2m}\times(0,\infty))$ to the equation
$$
\aligned&u_t=Lu+f(x,u)+\psi_{m,k}, \ (x,t)\in
B_{2m}\times(0,\infty);\\ &u(x,0)=g_m, \ x\in B_{2m};\\
&u(x,t)=0,\ (x,t)\in \partial B_{2m}\times(0,\infty),
\endaligned\tag4.2
$$
where $g_m\ge0$ is continuous and satisfies
$$
g_m(x)=\cases0,\text{ for}\ \ x\in B_m\\
m^2W,\ \text{for}\ x\in B_{2m}-B_{m+1}\endcases.
$$
(This construction is similar to the one in \cite{3,equation (2.5)}.)
Also,
$$
U(x,t)\equiv\lim_{m\to\infty}\lim_{k\to\infty}U_{m,k}(x,t)
\ \text{is a solution to (1.1) with initial condition}\ g=0.\tag4.3
$$
(See the two paragraphs after equation (2.6)  in \cite{3},
ignoring equation (2.6)
and the concept of a maximal solution that appears there.)

Consider (4.1)  with
$m$ replaced by $2m$, with the nonlinearity $f$ replaced by
$G$ as in condition (F-4a),
and with $g=W$.
Denote the solution to this equation by $u_m$.
We will show below that
$$
W-u_m\le U_{m,k}, \ \text{in}\ B_{\frac{3m}2}\times[0,\infty),
\text{for}\ k\ \text{sufficiently large, depending on }\ m.
\tag4.4
$$
Let $v_\infty$ denote the solution to $v'=G(v)$ with $v_\infty(0)=\infty$,
as in Lemma 1-i (note that by condition (F-4a), $G$ satisfies
the requirement in Lemma 1-i).
Since $f(x,u)\le G(u)$, we have
$f(x,v_\infty)-v_\infty'\le G(v_\infty)-v_\infty'=0$.
Also,
since $G$ is Lipschitz, it follows from the uniqueness theorem
for ordinary differential equations that $v_\infty(t)>0$,
for all $t\ge0$. Using these facts along with the
fact that $u_m=0$ on $\partial B_{2m}$ and the fact that
$v_\infty(0)=\infty$,
it follows from the maximum principle in Proposition 1
that
$$
u_m(x,t)\le v_\infty(t) \ \text{in}\  B_{2m}\times
(0,\infty).\tag4.5
$$
Letting $k\to\infty$ and then letting $m\to\infty$,
it follows from (4.3) that
the right hand side of (4.4) converges to
a solution $U$ of (1.1) with initial data $g=0$. By the uniqueness
assumption, $U=0$. Using this with (4.5)
then gives
$$
W(x)\le v_\infty(t) \ \text{in} \ R^n\times(0,\infty).\tag4.6
$$

We now  show that
$$\lim_{t\to\infty}v_\infty(t)=c_0,
\ \text{where}\ c_0\ \text{is the largest root of}\ G(u)=0.\tag4.7
$$
To see this, let $v_c$ be as in Lemma 1.
Integrating,  changing variables and letting $c\to\infty$, we obtain
$$
\int^\infty_{v_\infty(t)}\frac 1{-G(u)}du=t.\tag4.8
$$
Letting $t\to\infty$ in (4.8) and using the fact that $G$ is
locally Lipschitz
proves (4.7).
The theorem now follows from (4.6) and (4.7).

It remains to prove (4.4).
Let $V=W-u_m$. We have
$$
\aligned&LV+f(x,V)-V_t=(LW+f(x,W))-(Lu_m+G(u_m)-(u_m)_t)+\\
&f(x,W-u_m)-f(x,W)+G(u_m)
=\\
&f(x,W-u_m)-f(x,W)+G(u_m)\ge0\
\text{in}\ B_{2m}\times(0,\infty),
\endaligned
\tag4.9
$$
where the second equality follows from the definitions of $W$ and $u_m$,
and the  inequality follows from the definition of $G$.
On the other hand, we have
$$
LU_{m,k}+f(x,U_{m,k})-(U_{m,k})_t=-\psi_{m,k}\le0\ \text{in}\ B_{2m}\times(0,\infty).\tag4.10
$$
We now show that for sufficiently large $k$, depending on $m$,
$$
V(x,t)\le U_{m,k}(x,t), \ \text{on}\   \partial B_{\frac{3m}2}\times[0,\infty).
\tag4.11
$$
Define
$Q(x)=\left(l^2-(m+1+l-|x|)^2\right)W(x)$,
where $l=\frac12(m-1)$. Note that
$Q>0$ in the annulus $A_{m+1,2m}\equiv
\{m+1< |x|< 2m\}$
and vanishes on $\partial A_{m+1,2m}$.
Clearly
$LQ+f(x,Q)$ is bounded
in $A_{m+1,2m}\times[0,\infty)$.
Thus for $k$ sufficiently large,
we have
$LQ+f(x,Q)\ge-\psi_{m,k}$ in $A_{m+1,2m}\times
[0,T]$.
Since $Q(x)\le l^2W(x)<m^2W(x)=g_m(x)=U_{m,k}(x,0)$
on $A_{m+1,2m}$, and since
$Q$ vanishes on
$\partial A_{m+1,2m}$, it follows by the maximum
principle in Proposition 1 that
$U_{m,k}\ge Q$ in $A_{m+1,2m}\times[0,\infty)$, for
$k$ sufficiently large.
Substituting $|x|=\frac{3m}2$ in $Q$,
we conclude that for $m\ge4$ and  sufficiently large $k$,
$U_{m,k}(x,t)\ge Q(x)=(l^2-\frac14)W(x)> W(x)$
on $\partial B_{\frac{3m}2}\times[0,\infty)$.
This proves (4.11) since
$V\le W$.
In light of (4.9)-(4.11) and the fact that $V(x,0)=0$,
(4.4) now follows from the maximum principle in Proposition 1.
\hfill $\square$

\medskip

\noindent\bf Acknowledgement.\rm\ I thank Henri Berestycki
for sending me  the preprint \cite{1}. The  results there
concerning
the  existence/nonexistence of stationary solutions
in the case of periodic operators
proved to be a catalyst for  this paper.

\end